\documentclass[oneside,12pt]{amsart}
\usepackage{amssymb, amscd}

\newtheorem{thm}{Theorem}[section]
\newtheorem{cor}[thm]{Corollary}
\newtheorem{lem}[thm]{Lemma}
\newtheorem{defn}[thm]{Definition}
\newtheorem{prop}[thm]{Proposition}

\newtheorem{conj}[thm]{Conjecture}

\newcommand{\del}[2]{{}}

\title{Compatible Actions and Cohomology of Crystallographic Groups}

\author{Alejandro Adem$^*$}
\address{Department of Mathematics, University of British Columbia, Vancouver BC V6T 1Z2, Canada}%
\email{adem@math.ubc.ca}%

\author{Jianquan Ge}%
\address{Department of Mathematical Sciences, 
Tsinghua University, Beijing 100084, China}%
\email{gejq04@mails.tsinghua.edu.cn}%

\author{Jianzhong Pan }%
\address{Institute of Mathematics, Academia Sinica, Beijing 100080, China}%
\email{pjz@amss.ac.cn}%

\author{Nansen Petrosyan}
\address{Department of Mathematics, Indiana University, Bloomington IN 47405, USA}%
\email{nanpetro@indiana.edu}%

\thanks{$^*$The first author was partially supported by the NSF and
by NSERC}%

\date{\today}

\begin{document}
\begin{abstract}
We compute the 
cohomology of crystallographic 
groups $\Gamma= \mathbb Z^n \rtimes \mathbb Z/p$
with holonomy of prime order
by establishing the collapse at $E_2$
of the spectral sequence associated to their defining
extension. As an application we compute the group of
gerbes associated to
many
six--dimensional toroidal orbifolds arising in string
theory.
\end{abstract}

\maketitle

\section{Introduction}

Given a finite group $G$ and an integral representation
$L$ for $G$ (i.e. a homomorphism
$G\to GL_n(\mathbb Z)$, where $L$ is the underlying
$\mathbb ZG$--module), we can define the 
\textsl{semi--direct product} 
$\Gamma = L\rtimes G$.
Calculating the cohomology of these groups is a
problem of instrinsic algebraic interest; indeed 
if the representation is \textsl{faithful} then
these groups can be thought of
as crystallographic groups (see \cite{charlap},
page 74).

From the geometric point of view,
the action on $L$ gives rise to a $G$--action on the
$n$--torus $X=\mathbb T^n$; this approach can be
used to derive important examples
of orbifolds, known as \textsl{toroidal orbifolds} 
(see \cite{AP}).
In the case when $n=6$ these are of particular interest
in string theory (see \cite{bd}, \cite{EK}).

Given the split group extension $0\to L\to \Gamma\to G\to 1$,
the basic problem which we address is that of providing
conditions which imply the collapse (without extension
problems)
of the associated Lyndon--Hochshild--Serre spectral 
sequence. The conditions
which we establish are representation--theoretic, namely
depending solely on the structure of the integral representation
$L$. This can be a difficult problem (see \cite{totaro} for further
background) and there are well--known examples where the
spectral sequence does not collapse.

Our approach is to systematically apply 
the methods used in \cite{AP}. 
The key idea is to construct a free resolution $F$
for the semidirect product $L\rtimes G$ such that the 
Lyndon--Hochschild--Serre spectral sequence of the
group extension collapses at $E_2$. This requires
a chain--level argument, more specifically the
construction
of a \textbf{compatible} $G$--action on a certain free 
resolution for the torsion--free abelian group $L$
(see \S 2 for details).
We concentrate on the case of $G=\mathbb Z/p$, a cyclic
group of prime order $p$, as the representation theory
is well--understood. Our main algebraic 
result is the following

\begin{thm}\label{Big}
Let $G= \mathbb Z/p$, where $p$ is any
prime. If $L$ is any finitely generated 
$\mathbb ZG$--lattice\footnote{A $\mathbb ZG$--lattice
is a $\mathbb ZG$--module which happens to be a
free abelian group.},
and
$\Gamma = L\rtimes G$ is the associated semi--direct product
group, then for each $k\ge 0$
$$H^k(\Gamma , \mathbb Z)\cong \bigoplus_{i+j=k}
H^i(G,\wedge^j(L^*))$$
where $\wedge^j (L^*)$ denotes the $j$-th exterior
power of the dual module $L^*=Hom(L,\mathbb Z)$.
\end{thm}

Expressed differently: these results imply a complete calculation
for the integral cohomology of crystallographic groups
$\mathbb Z^n\rtimes \mathbb Z/p$ where $p$ is prime.
These calculations can be made explicit.

The theorem has an interesting geometric application:

\begin{thm} Let $G= \mathbb Z/p$, where $p$ is any prime.
Suppose that $G$ acts on a space $X$ homotopy equivalent to 
$(\mathbb
S^1)^n$ with $X^G\ne \emptyset$, then for each $k\ge 0$
$$H^k(EG\times_G X,\mathbb Z)\cong \bigoplus_{i+j=k}H^i(G,H^j(X,\mathbb Z))
\cong H^k(\Gamma, \mathbb Z).$$
where $\Gamma = \pi_1(X)\rtimes G$.

\end{thm}

On the other hand, the explicit computation for 
\textbf{torsion--free}
crystallographic groups with holonomy of prime order was
carried out long ago by Charlap and V\'asquez (see \cite{CV},
page 556).
Combining the two results we obtain a complete calculation:

\begin{thm}\label{complete}
Let $\Gamma$ denote a crystallographic group with holonomy
of prime order $p$, expressed as an extension
$$1\to L\to\Gamma\to \mathbb Z/p\to 1$$
where $L$ is a free abelian group of finite rank.
\begin{enumerate}
\item If $\Gamma$ is torsion--free, then $L\cong N\oplus \mathbb Z$
(it splits off a trivial direct summand) and
$$H^k(\Gamma, \mathbb Z)\cong H^0(\mathbb Z/p, \wedge^k(L^*))
\oplus H^1(\mathbb Z/p, \wedge^{k-1}(N^*))$$
for $0\le k\le rk(L)$; $H^k(\Gamma, \mathbb Z)=0$ for $k> rk(L)$.
\item If $\Gamma$ is not torsion--free, then
$H^*(\Gamma, \mathbb Z)$ can be computed using Theorem \ref{Big}.
\end{enumerate}
\end{thm}

In this paper we also consider the situation for
the cyclic group of order four; some partial results
are obtained but a general collapse has not been
established. However, based on these and other computations
we conjecture that for $G$ any cyclic group, the 
spectral sequence associated to a semi--direct product
of the form $\mathbb Z^n\rtimes G$
must collapse at $E_2$.

In the last section we give an application of our methods
to calculations for six--dimensional toroidal orbifolds,
showing that among
the 18 inequivalent $N=1$ supersymmetric string theories
on symmetric orbifolds of $(2,2)$--type without discrete
background, only two of them cannot be analyzed using our
methods i.e. we cannot show the existence of compatible
actions for the associated modules.  
If $\mathcal X=[\mathbb T^6/G]$  is
an orbifold arising this way, then our results provide 
a complete calculation for its associated group
of gerbes 
$Gb(\mathcal X)\cong H^3(EG\times_G \mathbb T^6,\mathbb Z)$
(see
\cite{AP} for more details).

\section{Preliminary Results}

The notion of a \textsl{compatible action} was first introduced
in \cite{br}. If such an action exists it allows one to
construct practical projective resolutions and from these to
compute the cohomology of the group. We will give the basic
definition and the main theorem that follows. More details can be
found in \cite{AP}.

Let $\Gamma = L \times_{\rho}G= L\rtimes G$ be the semidirect 
product of a
finite group $G$ and a finite dimensional $\mathbb Z$-lattice $L$
via a representation $\rho : G \rightarrow GL(L)$. $G$
acts on the group $L$ by the homomorphism $\rho$, and this extends
linearly to an action on the group algebra $R[L]$, where
$R$ denotes a commutative ring with unit. We write $l^g$
for $\rho (g)l$ where $l\in R[L], g\in G$.
In the rest of this paper $R$ will
represent $\mathbb Z$ (the integers) or  
$\mathbb Z_{(p)}$ (the ring of integers localized
at a fixed prime $p$).

\begin{defn}\label{T:comp}
Given a free resolution $\epsilon:F \to {R}$ of $R$ over $R[L]$,
we say that it admits an action of $G$ compatible with $\rho$ if
for all $g\in G$ there is an augmentation-preserving chain map
$\tau(g):F \to F$ such that

\begin{enumerate}
 \item{$\tau(g)[l\cdot f]=l^g\cdot [\tau(g)f]$
for all $l\in R[L]$  and $f\in F$},
\item{$\tau(g)\tau(g')=\tau(gg')$ for all $g,g'\in G$},
\item{$\tau (1) = 1_F$.}
\end{enumerate}
\end{defn}

The following two lemmas (see \cite{AP})
reduce the construction of
compatible actions to the case of faithful indecomposable
representations.

\begin{lem}\label{sum}
If $\epsilon_i: F_i \rightarrow  R$ is a projective $
R[L_i]$-resolution of $R$ for $i=1,2$, then $\epsilon_1 \otimes
\epsilon_2: F_1 \otimes F_2 \rightarrow R$ is a projective $ R[L_1
\times L_2]$-resolution of $R$. Furthermore, if $G$ acts
compatibly on $F_i$ by $\tau_i$ for $i=1,2$, then a compatible
action of $G$ on $\epsilon_1 \otimes \epsilon_2: F_1 \otimes F_2
\rightarrow R$ is given by $\tau(g)(f_1\otimes f_2)=\tau_1(g)(f_1)
\otimes \tau_2(g)(f_2)$.
\end{lem}

\begin{lem}
If $L$ is a $ R[G_1]$-module, $\pi: G_2 \rightarrow G_1$ a group
homomorphism, and $\epsilon:F\rightarrow R$ is a $
R[L]$-resolution of $R$ such that $G_1$ acts compatibly on it by
$\tau^{\prime}$, then $G_2$ also acts compatibly on it by
$\tau(g)f=\tau^{\prime}(\pi(g))f$ for any $g\in G$.
\end{lem}

If a compatible action exists, we can give $F$ a $\Gamma$-module
structure as follows. An element
$\gamma\in \Gamma$ can be
expressed uniquely as $\gamma=lg$, with $l\in L$ and $g\in G$. We
set $\gamma\cdot f=(lg)\cdot f=l\cdot \tau(g)f$. Note that given
any $G$--module $M$, this inflates to a $\Gamma$--action on $M$
via the projection $\Gamma\to G$.

We can always construct a \textit{special} free
resolution $F$ of $R$ over $L$, characterized by the property that
the cochain complex $Hom_L(F,R)$ for computing the cohomology 
$H^\ast(L,R)$
has all coboundary maps zero (more details will be provided
in the next section). 
Using this fact, the following
was proved in 
\cite{AP}:

\begin{thm}\label{T:2.3}{\normalfont (Adem-Pan)}
Let $\epsilon:F \to{R}$ be a special free resolution of ${R}$ over
$L$ and suppose that there is a compatible action of $G$ on $F$.
Then for all integers $k\ge 0$, we have
$$H^k(L \rtimes G, R)=\bigoplus_{i+j=k}H^i(G,H^j(L, R)).$$
\end{thm}

This result can be interpreted as saying that the 
Lyndon-Hochschild-Serre spectral sequence
$$E^{p,q}_2=H^p(G,H^q(L, R))\Rightarrow H^{p+q}(L\rtimes G, R)$$
collapses at $E_2$ without extension problems. Note that
this is not always the case; in fact
there are examples of semi--direct products of the form
$\mathbb Z^n\rtimes (\mathbb Z/p)^2$ where the associated
spectral sequence has non--trivial differentials (see \cite{totaro}).
This will be discussed in \S 5.

\section{Construction of Compatible Actions}

Let $R[L]$ denote the 
group ring of $L$, a free abelian group with
basis $\{x_1,\dots ,x_n\}$. Then the elements 
$x_1-1, \dots , x_n-1$ form a \textsl{regular sequence}
in $R[L]$, hence the Koszul complex $K_*=K(x_1-1,\dots, x_n-1)$
is a free resolution of the trivial module $R$. It has
the additional property of being a \textsl{differential
graded algebra} (or DGA). We briefly recall how it looks.
There are generators $a_1, \dots, a_n$ in degree one,
and the graded basis for $K_*$ can be identified with
the usual basis for the exterior algebra they generate.
The differential is given
by the following formula: if $a_{i_1\dots i_p}=a_{i_1}\dots a_{i_p}$
is a basis element in $K_p$, then

$$d(a_{i_1\dots i_p}) = 
\sum_{j=1}^p (-1)^{j-1}(x_{i_j}-1)a_{i_1\dots\hat{i_j}\dots i_p}.$$

Now the cohomology of the free abelian group $L$ is 
precisely an exterior algebra on $n$ one-dimensional
generators, which in fact can be identified with the
dual elements $a_1^*,\dots ,a_n^*$. In particular
we see that the cochain complex $Hom_{R[L]}(K_*,R)$
has zero differentials, and hence $K_*$ is a 
\textbf{special} free resolution of $R$ over $R[L]$
(this resolution also appears in \cite{charlap} pp. 
96--97).

We now consider how to construct a compatible
$G$--action on $K_*$, given a $G$--module structure
on $L$.

\begin{thm} If $G$ acts on the
lattice $L$, let $K_*= K(x_1-1,\dots x_n-1)$ 
denote the special
free resolution of $R$ over $R[L]$ defined using the Koszul
complex associated to the elements $x_1-1,\dots ,x_n-1$,
where $\{x_1,\dots ,x_n\}$ form a basis for $L$.
Suppose that there is a homomorphism $\tau : G\to Aut(K_1)$ such
that for every $g\in G$ and $a\in K_1$ 
it satisfies
$$ d\tau(g)(a)= d(a)^g$$
where $d:K_1\to K_0$ is the usual Koszul differential, and
$d(a)^g\in K_0=R[L]$.
Then $\tau$ extends to $K_*$ using its DGA structure
and so defines a compatible $G$--action on $K_*$.
\end{thm}
\begin{proof}
First we observe that $\tau(g)$ acts on $K_0=R[L]$ via the
original $G$--action, i.e. $\tau(g) (x) = x^g$ for any
$x\in K_0$. Next we define the action on the basis of
$K_*$ as a graded $R[L]$--module, namely:
$$\tau (g) (a_{i_1}\dots a_{i_p}) 
=\tau(g)(a_{i_1})\dots\tau (g)(a_{i_p}).$$
If $\alpha\in R[L]$ and $u\in K_*$, we
define
$\tau (g) (\alpha u) = \alpha^g\tau(g) (x)$.
By linearity and the DGA structure of $K_*$, 
this will define $\tau: G \to Aut(K_*)$, with the
desired properties.
\end{proof}

%
%

Generally speaking it can be quite difficult to construct a
compatible action; however there is an important special
case where it is quite straightforward.

\begin{thm}\label{RG}
Let $\phi: G\to \Sigma_n$ denote a group homomorphism, 
where $\Sigma_n$ denotes the symmetric
group on $n$ elements. Let $G$ act on $\mathbb Z^n$ 
via\footnote{Such a module is called a \textsl{permutation
module}.} this
homomorphism. Then the associated Koszul complex $K_*$
admits a compatible $G$--action.
\end{thm}
\begin{proof}
By Lemma 2.3 we can assume that $G$ is a subgroup of $\Sigma_n$,
hence it will suffice to prove this for
$\Sigma_n$ itself.
If we take generators $a_1,\dots , a_n$ for the Koszul complex
corresponding to the elements $x_1, \dots , x_n$ in the underlying
module $L$, then we can define $\tau$ as follows: 
$$\tau (\sigma ) (a_i) = a_{\sigma (i)}.$$
This obviously defines a permutation representation on
$K_1$, and compatibility follows from the fact that
for all $a_i$, $1\le i\le n$ and $\sigma\in\Sigma$ we have

$$d\tau (\sigma)(a_i)=d(a_{\sigma (i)})=x_{\sigma (i)}-1
=(x_i-1)^{\sigma}.$$
This completes the proof.
\end{proof}

Aside from permutation representations, it is difficult to
construct general examples of compatible actions. However
if $G$ is a cyclic group then we can handle an important
additional type of module.

\begin{prop}\label{IG}
Let the cyclic group $G=\left< t|t^n=1 \right>$ act on $\mathbb
Z^{n-1}$ by:

\[\xi_1: t\mapsto \left( \begin{array}{llllll} 0 &1 & 0 & \dots & 0 & 0 \\
0 & 0 & 1 & \dots & 0 & 0\\  &  & \dots & &\\0 & 0 & 0& \dots & 1 & 0
\\-1 & -1 & -1 & \dots & -1 & -1\end{array}\right)\in GL_{n-1}(\mathbb Z).
\]
If $x_1,\dots ,x_n$ is the canonical basis under which
the action is represented by the matrix above,
then the free resolution $K_*=K_*(x_1-1,\dots , x_n-1)$ 
admits an action of
$G$ compatible with $\xi_1$, which can be defined by:
\[ \tau(t)(a_1)= -x_{n-1}^{-1}a_{n-1} 
\qquad \tau(t)(a_k)=-x_{n-1}^{-1}(a_{n-1}-a_{k-1}), ~~~1< k \leq
n-1.\]
\end{prop}

\begin{proof}
The proof is a straightforward calculation 
verifying that 
$\tau$ defines a compatible action. First we verify
that $\tau^n = 1$. For this we observe that if $A$
is the matrix in $GL_{n-1}(\mathbb Z)$ representing
the generator $t$, then expressed in terms of the
basis $\{a_1,\dots ,a_n\}$ we have that 
$\tau (t) = x_{n-1}^{-1} A$. If we iterate this
action and use the fact that $\tau(g) (\alpha u)
= \alpha^g \tau (g)(u)$ then we obtain

$$\tau (t)^n = (x_{n-1}^{-1})^{t^{n-1}}(x_{n-1}^{-1})^{t^{n-2}}
\dots (x_{n-1}^{-1})^t(x_{n-1}^{-1})A^n = 1.$$
This follows from the fact that the characteristic
polynomial of $A$ is the cyclotomic polynomial
$p(z)=1 + z + \dots + z^{n-1}$, hence we have
that $p(A)=0$ on the underlying module $L$ 
and so in multiplicative notation
we have that $u^{t^{n-1}}\cdot u^{t^{n-2}}\cdot\dots\cdot u^t\cdot u=1$
for any $u\in L$.

Next we verify compatibility:
$$\tau (t) d(a_1) = \tau(t)(x_1-1)=x_{n-1}^{-1}-1$$
$$d\tau(t)(a_1) = d(-x_{n-1}^{-1}a_{n-1})= -x_{n-1}^{-1}
(x_{n-1}-1)=x_{n-1}^{-1}-1.$$
Similarly 
for all $1<
k\leq n-1$ we have that:

\[\tau(t)
d(a_k)=\tau (t)(x_k-1)=x_{k-1}x_{n-1}^{-1}-1\]
\[=-x_{n-1}^{-1}(x_{n-1}-1-x_{k-1}-1)=d(-x_{n-1}^{-1}(a_{n-1}-a_{k-1}))
=d\tau(t)(a_k).\]

\end{proof}

%
%
%

For $G=\mathbb Z/n$,
the module which gives rise to the matrix in
\ref{IG} is the augmentation ideal $IG$, which has
rank equal to $n-1$. 
The following proposition is an application of the 
results in this section.

\begin{prop}\label{(r,s,t)}

Let $G=\mathbb Z/n$, and
assume that $L$ is a $\mathbb ZG$--lattice
such that
$$L\cong M
\oplus IG^t$$ 
where $M$ is a permutation module.
Then, for any coefficient ring $R$, 
the special free resolution
$K_*$ over $R[L]$ admits a compatible
$G$--action.
\end{prop}

\begin{proof}
This follows from applying Lemma \ref{sum} to 
Theorem \ref{RG} and Proposition \ref{IG}.
\end{proof}

For our cohomology calculations it will be practical to
use the coefficient ring $R=\mathbb Z_{(p)}$, where
$p$ is a prime.
In this situation, for $G=\mathbb Z/p$ 
(see \cite{CR}) 
there are only three distinct 
isomorphism classes of indecomposable
$RG$--lattices, namely $R$ (the trivial module), $IG$ (the
augmentation ideal) and $RG$, the group ring. 
Moreover, if $L$ is \textsl{any} finitely generated
$\mathbb ZG$--lattice, 
we can construct
a $\mathbb ZG$-homomorphism 
$f:L'\to L$ such that
\begin{itemize}
\item $L'\cong \mathbb Z^r\oplus \mathbb ZG^s\oplus IG^t$
\item $f$ is an isomorphism after tensoring with $R$.
\end{itemize}
We shall call $L'$ a representation of type $(r,s,t)$.

\section{Applications to Cohomology}

We are now ready to prove our main result.

\begin{thm}\label{mainthm}Let $G= \mathbb Z/p$, where $p$ is any
prime. If $L$ is any finitely generated $\mathbb ZG$--lattice,
and 
$\Gamma = L\rtimes G$ is the associated semi--direct product
group, then for each $k\ge 0$
$$H^k(\Gamma , \mathbb Z)\cong \bigoplus_{i+j=k}
H^i(G,\wedge^j(L^*))$$
where $\wedge^j (L^*)$ denotes the $j$-th exterior
power of the dual module $L^*=Hom(L,\mathbb Z)$.
\end{thm}

\begin{proof}

First, let us prove the analogous result for the
cohomology with coefficients in $R=\mathbb Z_{(p)}$. 
We make the assumption that $L$ is a module of
type $(r,s,t)$.
We need to verify:
\[H^k(\Gamma , R)\cong \bigoplus_{i+j=k} H^i(G,\wedge^j(L_R^*))
\]
where $L_R^* = L^*\otimes R$.
In fact, we see that in the associated
Lyndon-Hochschild-Serre spectral sequence for the extension
$0\to L\to \Gamma\to G\to 1$ with 

\[
E^{i,j}_2(R) =H^i(G,H^j(L, R))
\Rightarrow H^k(\Gamma, R)
\]
there are no differentials and no extension problems.
This follows from applying Theorem \ref{T:2.3}
and the fact that
the module $L$ gives rise to a special resolution with a
compatible action by Proposition \ref{(r,s,t)}.

Now let us consider the case when $L$ is not of type
$(r,s,t)$. As observed previously, we
can construct a $\mathbb ZG$--lattice
$L'$ and a map $f:L'\to L$ such that 
$L'$ is of type $(r,s,t)$
and $f$ is an isomorphism after tensoring with $R$.
Under these conditions
$f$ will induce a map between the spectral sequences
with $R$--coefficients
for the extensions corresponding to $L$ and $L'$.
However by our hypotheses, $\wedge^k(L'^*_R)$
and $\wedge^k(L^*_R)$ are isomorphic 
as $RG$--modules for all $k\ge 0$,
with the isomorphism induced by $f$. Hence the
corresponding $E_2$--terms are isomorphic,
and so the spectral sequences both collapse
and the result follows.

It now remains to prove the result with coefficients in
the integers $\mathbb Z$. Note that by the universal
coefficient theorem, we have
$H^*(\Gamma, \mathbb Z_{(p)})
\cong H^*(\Gamma, \mathbb Z)\otimes\mathbb Z_{(p)}$
hence the only relevant 
discrepancy between $H^*(\Gamma, \mathbb Z_{(p)})$
and $H^*(\Gamma, \mathbb Z)$ might arise from the presence
of torsion prime to $p$
in the integral cohomology of $\Gamma$.
However, a quick inspection of the 
spectral sequence of
the extension $0\to L\to \Gamma\to G\to 1$ with
$\mathbb Z_{(q)}$ coefficients shows that
there is no torsion prime to $p$ in the cohomology, as
$L$ is free abelian and $G$ is a $p$--group. This
completes our proof.
\end{proof}

In \cite{AP}, Corollary 3.3, it was observed that the
spectral sequence for the extension $L\rtimes G$ satisfies
a collapse at $E_2$ wthout extension problems 
if the same is true for all the
restricted extensions $L\rtimes G_p$, where the $G_p\subset G$
are the $p$--Sylow subgroups of $G$. We obtain the following

\begin{cor}\label{composite}
Let $G$ denote a finite group of square--free order,
and $L$ any finitely generated $\mathbb Z G$--lattice.
Then for all $k\ge 0$ we have

$$H^k(L\rtimes G, \mathbb Z)
\cong \bigoplus_{i+j=k} H^i(G, \wedge^j(L^*)).$$
\end{cor}

We now consider a more geometric situation. Suppose that the
group $G=\mathbb Z/p$ acts on a space $X$ which has the 
homotopy type of a product of circles.

\begin{thm} Let $G= \mathbb Z/p$, where $p$ is any prime.
Suppose $G$ acts on a space $X$ homotopy equivalent to $(\mathbb
S^1)^n$ with $X^G\ne \emptyset$, then for each $k\ge 0$
$$H^k(EG\times_G X,\mathbb Z)\cong \bigoplus_{i+j=k}H^i(G,H^j(X,\mathbb Z))
\cong H^k(\Gamma, \mathbb Z)$$
where $\Gamma = \pi_1(X)\rtimes G$.
\end{thm}

\begin{proof} The space $EG\times_G X$ fits into a fibration
$X\hookrightarrow EG\times_G X \to BG$
which has a section due to the fact that 
$X^G\ne \emptyset$. 
Let $\Gamma$ denote the fundamental group of $EG\times_G X$. 
The long
exact sequence for the homotopy groups of the fibration gives
rise to a split extension
$1\to \pi_1(X)\to
\Gamma \to
G\to 1$.
Since $\pi_1(X)\cong
L$, a $\mathbb ZG$--lattice, 
this shows that $\Gamma\cong L\rtimes G$, where the $G$ action
is induced on $L$ via the action on the fiber.
Note that $EG\times_G X$ is an Eilenberg-MacLane space of type
$K(\Gamma, 1)$. Hence, 
$H^*(EG\times_G X, \mathbb Z)\cong
H^*(\Gamma, \mathbb Z)$ 
and the result follows from 
Theorem \ref{mainthm}.

\end{proof}

Note that a special case of this result was proved in
\cite{adem}, namely for actions where $\pi_1(X)\otimes\mathbb Z_{(p)}$ 
is isomorphic
to a direct sum of indecomposables of rank $p-1$.
The
terms $H^i(\mathbb Z/p, \wedge^j(L^*))$ can be computed
if $L$ is known up to isomorphism. In fact, all we need
is to know $L$ up to $\mathbb Z/p$ cohomology, as this
will determine its indecomposable factors (at least up
to $\mathbb Z_{(p)}$--equivalence).

As we mentioned in the introduction, our results complete the
calculation for the cohomology of crystallographic groups
with prime order holonomy when combined with previous work
on the torsion--free case (Bieberbach groups). The terms
appearing in the formulas in Theorem \ref{complete} can be
explicitly computed, as was observed in \cite{CV}.

\section{Extensions to Other Groups}

In this section we explore to what degree our results
can be extended to other groups. In the case
of the cyclic group of order four, the indecomposable
integral representations are easy to describe, so it
is a useful test case.

Let $G=\mathbb Z/4$, from \cite{BG} we can give a complete
list of all (nine) indecomposable pairwise nonequivalent integral
representations by the following adopted table, where $a$
is a generator for $\mathbb Z/4$:

$$
\rho_1:a\rightarrow 1; \rho_2:a\rightarrow -1;
\rho_3:a\rightarrow \left( \begin{array}{ll} 0 & 1
\\ 1 & 0 \end{array}\right); 
\rho_4:a \rightarrow \left( \begin{array}{ll} 0 & 1
\\ -1 & 0 \end{array}\right); $$

$$\rho_5:a \rightarrow \left( \begin{array}{lll} 0 & 0 & -1 \\
1 & 0 & -1\\0 & 1 & -1 \end{array}\right);\qquad
\rho_6: a\rightarrow \left( \begin{array}{lll} 0 & 1 & 0 \\
-1 & 0 & 1\\0 & 0 & 1 \end{array}\right)$$

$$\rho_7:a \rightarrow \left( \begin{array}{llll} 0 & 1& 0 & 0\\
 0 & 0 & 1 & 0\\0 & 0 & 0 & 1\\ 1 & 0 & 0 & 0 \end{array}
\right);\qquad
\rho_8:a \rightarrow \left( \begin{array}{llll} 0 & 0 & -1 & 1\\
1 & 0 & -1 & 1\\0 & 1 & -1 & 0\\ 0 & 0 & 0 & 1
\end{array}\right)$$

$$\rho_9:a \rightarrow \left( \begin{array}{llll} 0 & 1& 0 & 0\\
-1 & 0 & 0 & 1\\0 & 0 & -1 & 1\\ 0 & 0 & 0 & 1 \end{array}
\right)$$

\begin{thm}\label{Z/4} Let $G= \mathbb Z/4$ and $L$ a
finitely generated $\mathbb ZG$--lattice. 
If $L$ is a direct sum of indecomposables of type
$\rho_i$ for $i\leq 7$, and $i\ne 6$, then there is a
compatible action and 

\[H^k(L\rtimes G, \mathbb Z)=
\bigoplus_{i+j=k}H^i(G,H^j(L, \mathbb Z)).\]
\end{thm}

\begin{proof} 
For the indecomposables $\rho_1, \rho_2,\rho_3,\rho_4$, compatible
actions are known to exist on the associated resolutions by
the results in \cite{AP}. The same is true\footnote{In fact
$\rho_5$ corresponds to the dual module $IG^*$, but for cyclic
groups there is an isomorphism $IG\cong IG^*$.} for $\rho_5$ and
$\rho_7$ by \ref{IG} and \ref{RG}. Hence if $L$ is any integral
representation expressed as a direct sum of $\rho_i$,
$i\le 7$ and $i\ne 6$, the result follows from \ref{sum} and 
\ref{T:2.3}.
\end{proof}

In the case of $\rho_6$, $\rho_8$ and $\rho_9$, 
a compatible action is not
known to exist. However, via an explicit
computation done in \cite{petrosyan}, 
we can establish the collapse of the
spectral sequence for the extension $\Gamma_6$
associated to
$\rho_6$, yielding

\[H^i(\Gamma_6, \mathbb Z)= \left\{ \begin{array}{ll}
\mathbb{Z} & \mbox{ if $i=0,1$ }\\ \mathbb Z/4 \oplus \mathbb Z &
\mbox{ if $i=2$ } \\ \mathbb Z/2\oplus \mathbb Z & \mbox{ if $i=3$ }\\
\mathbb Z/4 \oplus \mathbb Z/2& \mbox{ if $i\geq 4$}
\end{array} \right.
 \]
which verifies the statement analogous to \ref{Z/4} for the
cohomology of $\Gamma_6$.

Indeed, for all the examples of semidirect products we 
have considered so
far, there is a collapse at $E_2$ in the Lyndon-Hochschild-Serre
spectral sequence of the group extension
$0\rightarrow L\rightarrow G\rtimes L \rightarrow G
\rightarrow 1$
and therefore we can make the following:

\begin{conj} Suppose that $G$ is a finite cyclic group and $L$ a 
finitely generated $\mathbb Z G$--lattice; then
for any $k\geq 0$ we have
$$H^k(G\rtimes L, \mathbb Z)=
\bigoplus_{i+j=k}H^i(G,H^j(L, \mathbb Z)).$$
\end{conj}

In \cite{totaro} examples were given of semi-direct products
of the form $L\rtimes (\mathbb Z/p)^2$ where the associated
mod $p$ Lyndon-Hochschild-Serre
spectral sequence has non--zero differentials. This relies on
the fact that for $G=\mathbb Z/p\times\mathbb Z/p$, there
exist $\mathbb ZG$--modules $M$ which are not realizable as
the cohomology of a $G$--space. These are the counterexamples
to the \textsl{Steenrod Problem} given by G.Carlsson
(see \cite{carlsson}), where $M$ can be identified with $L^*$.
No such counterexamples exist for finite cyclic groups
$\mathbb Z/N$, which means that disproving our conjecture
will require a different approach.

We should also mention that by using results due 
to Nakaoka (see \cite{evens}, pages 19 and 50)
we know that the spectral sequence for a wreath product 
$\mathbb Z^n\rtimes G$, where $G\subset\Sigma_n$ (the symmetric
group) acting on
$\mathbb Z^n$ via permutations will always
collapse at $E_2$, without extension problems. This can be
interpreted as the fact that a strong collapse theorem 
holds
for all \textsl{permutation modules} and all finite groups $G$.
A simple proof of this result can be obtained by applying
Theorem \ref{T:2.3} to Proposition \ref{RG}.

As suggested by
\cite{totaro}, the results here can be considered
part of a very general problem, which is both
interesting and quite challenging:

\bigskip

\noindent\textbf{Problem}: \textsl{Given a finite group $G$,
find suitable conditions on
a $\mathbb ZG$--lattice $L$ so that the spectral sequence
for $L\rtimes G$ collapses at $E_2$}.

\section{Application to Computations for Toroidal Orbifolds}

Interesting examples
arise from calculations for six--dimensional orbifolds, where the
usual spectral sequence techniques become rather complicated. Here our
methods provide an important new ingredient that allows us to compute
rigorously beyond the known range. 
An important class of examples in physics arises
from actions of a cyclic group
$G=\mathbb Z/N$ on $\mathbb T^6$.
In our scheme, these come from
six-dimensional integral representations of $\mathbb Z/N$. However,
the constraints
from physics impose certain restrictions on them (see \cite{EK}, 
\cite{vw}).
If $\theta\in GL_6(\mathbb Z)$ is an element of order $N$, then
it can be diagonalized over the complex numbers.
The associated eigenvalues, denoted
$\alpha_1,\alpha_2,\alpha_3$, should satisfy
$\alpha_1\alpha_2\alpha_3 =1$, and in addition all of the
$\alpha_i\ne 1$. The first condition implies that
the orbifold $\mathbb T^6\to \mathbb T^6/G$ is a
\textsl{Calabi--Yau orbifold}, and so admits
a \textsl{crepant resolution}.
These more restricted representations
have been classified\footnote{In the language of
physics, they show that there exist 18 inequivalent $N=1$ supersymmetric
string theories on symmetric orbifolds of $(2,2)$--type without
discrete background.} in \cite{EK},
where it is shown that there are precisely 18 inequivalent lattices
of this type.

It turns out that calculations are focused on computing 
the equivariant
cohomology 
$H^*(EG\times_G\mathbb T^6, \mathbb Z)$ (see \cite{AP} and \cite{bd}
for more details). 
As was observed in \cite{AP}, we can compute the group
of gerbes associated to the orbifold
$\mathcal X = [\mathbb T^6/G]$ via the
isomorphism
$$Gb(\mathcal X) \cong H^3(EG\times_G\mathbb T^6, \mathbb Z)
\cong H^3(\mathbb Z^6\rtimes G,\mathbb Z),$$
whence our methods can be used to
obtain some fairly complete results in this setting.
Before proceeding we recall that as in Corollary \ref{composite}
the collapse of the spectral sequence
for an extension $L\rtimes G$ will follow from the
existence of compatible $Syl_p(G)$ actions on the
Koszul complex $K_*$ for every
prime $p$ dividing $|G|$. If these exist we shall say
that $K_*$ admits a \textsl{local} compatible action.
\begin{thm}
Among the $18$ inequivalent integral representations associated
to the six--dimensional orbifolds $\mathbb T^6/\mathbb Z/N$
described above, only two of
them are not known to admit (local) compatible actions. Hence for those
$16$ examples there is an isomorphism

$$H^k(E\mathbb Z/N\times_{\mathbb Z/N}\mathbb T^6, \mathbb Z)
\cong \bigoplus_{i+j=k} 
H^i(\mathbb Z/N, H^j(\mathbb T^6,\mathbb Z))$$
\end{thm}

\begin{proof}

Consider the defining matrix of an indecomposable
action of $\mathbb{Z}/N$ on $\mathbb{Z}^n$ with determinant
one, expressed in 
canonical form
as

$$\theta=\left(\begin{array}{cccccc}
               0&&&\ldots&0&v_1\\
               1&0&&\ldots&0&v_2\\
               0&1&0&\ldots&0&v_3\\
               \vdots&&\ddots&&&\vdots\\
               0&&\ldots&&1&v_n\end{array}\right)$$
where $v_1=\pm 1$.

In \cite{EK} it was determined that the matrices that
specify the indecomposable modules appearing as
summands for the $N=1$
supersymmetric $\mathbb Z/N$-orbifolds can be given as follows,
where the vectors represent the values $(v_1,v_2,\dots , v_n)$:

\vskip  .3 cm \begin{center}{\textbf{Indecomposable matrices relevant 
for 
$N=1$ supersymmetry}}\end{center}
 \vskip .3 cm
\begin{center}
\begin{tabular}{|l|l|l|l|}
\hline \multicolumn{1}{|c|}{$n=1$} &\multicolumn{1}{|c|}{$n=2$}
&\multicolumn{1}{|c|}{$n=3$} &\multicolumn{1}{|c|}{$n=4$}
 \\[.31 mm]
\hline
\hline &&&\\[-3.5 mm]
$\mathbb Z/2^{(1)}:\,\,(-1)$& $\mathbb Z/3^{(2)}:\,\,(-1,-1)$&
$\mathbb Z/4^{(3)}:\,\,(-1,-1,-1)$& $\mathbb Z/6^{(4)}:\,\,(-1,0,-1,0)$
\\[ 2 mm]
& $\mathbb Z/4^{(2)}:\,\,(-1,0)$& $\mathbb Z/6^{(3)}:\,\,(-1,0,0)$&
$\mathbb Z/8^{(4)}:\,\,(-1,0,0,0)$
\\[ 2 mm]
& $\mathbb Z/6^{(2)}:\,\,(-1,1)$& & $\mathbb Z/{12}^{(4)}:\,\,(-1,0,1,0)$
\\
\hline
\end{tabular}
\vskip  .5 cm
\begin{tabular}{|l|l|}
\hline \multicolumn{1}{|c|}{$n=5$} &\multicolumn{1}{|c|}{$n=6$}
 \\[.31 mm]
\hline
\hline &\\[-3.5 mm]
$\mathbb Z/6^{(5)}:\,\,(-1,-1,-1,-1,-1)$&
$\mathbb Z/7^{(6)}:\,\,(-1,-1,-1,-1,-1,-1)$
\\[2 mm]
$\mathbb Z/8^{(5)}:\,\,(-1,-1,0,0,-1)$& $\mathbb Z/8^{(6)}:\,\,(-1,0,-1,0,-1,0)$
\\[2 mm]
& $\mathbb Z/{12}^{(6)}:\,\,(-1,-1,0,1,0,-1)$
\\[2 mm]
\hline
\end{tabular}
\end{center}

We will show that all of these, except possibly
$\mathbb Z/8^{(5)}$ and $\mathbb Z/{12}^{(6)}$, admit 
local compatible actions.
The examples of rank two or less were dealt with in
\cite{AP}; for $N =2, 3, 6, 7$ the result follows directly
from \ref{mainthm} and \ref{composite}. 
The case $\mathbb Z/4^{(3)}$ was covered in \ref{Z/4}.
We will deal explicitly with the cases
$\mathbb Z/8^{(4)}$, $\mathbb Z/{12}^{(4)}$ and $\mathbb Z/8^{(6)}$.
\medskip

\noindent (1)~ The group
$\mathbb Z/8$ acts on $\mathbb Z^4$ with generator represented
by the matrix:
$$T=\left(\begin{array}{cccc}0& 0& 0& -1\\1& 0& 0& 0 \\ 0& 1&
0& 0\\0& 0& 1& 0\end{array}\right)$$

\noindent We define a compatible action by the following formulas:
$$\tau(t)(a_1)=-x_4^{-1}a_4,
\tau(t)(a_2)=a_1,
\tau(t)(a_3)=a_2,
\tau(t)(a_4)=a_3.$$

\noindent (2)~ The group
$\mathbb Z/{12}$ acts on $\mathbb Z^4$ with generator represented
by the matrix:
$$T=\left(\begin{array}{cccc}0& 0& 0& -1\\1& 0& 0& 0 \\ 0& 1&
0& 1\\0& 0& 1& 0\end{array}\right)$$
For this example
it suffices to construct a compatible action for $\mathbb Z/4$ (with
generator represented by the matrix $T^3$) as we already
know that a compatible action exists restricted to $\mathbb Z/3$.
Now we have that $\mathbb Z/4$ acts on $\mathbb Z^4$ with
$$T^3=\left(\begin{array}{cccc}0& -1& 0& -1\\0& 0& -1& 0\\0& 1& 0& 0\\
1& 0& 1& 0\end{array}\right)$$
which is a matrix whose square is $-I$. This implies
that the module is a sum of two copies of the faithful
rank two
indecomposable (see \S 5), for which a compatible action is known
to exist (as explained in Theorem \ref{Z/4}), and so this case
is taken care of.
\medskip

\noindent (3)~The group $\mathbb Z/8$ acts on $\mathbb Z^6$ with 
generator represented
by the matrix:
$$T=\left(\begin{array}{cccccc}0& 0& 0& 0& 0& -1\\1& 0& 0& 0& 0& 0\\
0& 1& 0& 0& 0& -1\\0& 0& 1& 0& 0& 0\\0& 0& 0& 1& 0& -1\\0& 0& 0& 0& 1& 0
\end{array}\right)$$
The formulas for a compatible action
are given by
$$\tau(t)(a_1)=-x_6^{-1}a_6,
\tau(t)(a_2)=a_1,
\tau(t)(a_3)=x_6^{-1}(a_2-a_6)$$
$$\tau(t)(a_4)=a_3,
\tau(t)(a_5)=x_6^{-1}(a_4-a_6),
\tau(t)(a_6)=a_5.$$
We have shown that (local) 
compatible actions exist for all representations
constructed using indecomposables other than $\mathbb Z/8^{(5)}$
and $\mathbb Z/12^{(6)}$. However, these indecomposables can only
appear once in the list due to dimensional constraints, namely in
the form $\mathbb Z/8^{(5)}\oplus \mathbb Z/2^{(1)}$ and
$\mathbb Z/12^{(6)}$ itself.
Thus our proof is 
complete.
\end{proof}

\end{document}